\documentclass[12pt]{article}
\usepackage{amssymb,amsmath,latexsym,amsthm, color}
\usepackage{amssymb,amsmath,latexsym}
\usepackage{color,amsmath,amsfonts,amssymb,epsfig,latexsym,pstricks,xkeyval, mathabx}
\usepackage[utf8]{inputenc}
\usepackage[T1]{fontenc}
\usepackage[english]{babel}
\usepackage{esint}
\usepackage{amssymb, amsthm, bbm, amsmath, latexsym, color}
\usepackage{enumerate, mathrsfs, hyphenat, graphicx}
\usepackage[font=sf, labelfont=sf]{caption}
\usepackage{thmtools, thm-restate}
\usepackage{thm-patch, thm-kv, thm-autoref}
\usepackage[colorlinks=true, urlcolor=blue, linkcolor=blue, citecolor=black]{hyperref}
\theoremstyle{plain}
\declaretheorem[title=Theorem, parent=section]{theorem}
\declaretheorem[title=Lemma,sibling=theorem]{lemma}

\declaretheorem[title=Proposition,sibling=theorem]{proposition}
\declaretheorem[title=Remark,sibling=theorem]{remark}

\theoremstyle{definition}

\numberwithin{equation}{section}
\addto\extrasenglish{%

}

\usepackage[backgroundcolor=none, bordercolor=blue]{todonotes}
\makeatletter
\providecommand\@dotsep{5}
\renewcommand{\listoftodos}[1][\@todonotes@todolistname]{%
  \@starttoc{tdo}{#1}}
\makeatother

\allowdisplaybreaks

\begin{document}
\author{
 {\bf Zoran Vondra\v{c}ek}
\thanks{Research supported in part by the Croatian Science Foundation under the project 4197.}
}

\title{\bf A probabilistic approach to a non-local quadratic form and its connection to the Neumann boundary condition problem}

\date{}

\maketitle

\begin{abstract}
In this paper we look at a probabilistic approach to a non-local quadratic form that has lately attracted some interest. This form is related to a recently introduced non-local normal derivative. The goal is to construct two Markov process: one corresponding to that form and the other which is related to a probabilistic interpretation of the Neuman problem. We also study the Dirichlet-to-Neumann operator for non-local operators. 
\end{abstract}



\section{Introduction}\label{s:intro}
Let $D\subset {\mathbb R}^d$, $d\ge 2$, be a bounded open set. For $\alpha\in (0,2)$ set $k(x,y)=c(d,\alpha)|x-y|^{-d-\alpha}$, $x,y\in {\mathbb R}^d$, where $c(d,\alpha)$ is a positive constant. Consider the symmetric bilinear form
\begin{equation}\label{e:form-hat}
\widehat{{\mathcal E}}(u,v):=  \frac12 \iint_{({\mathbb R}^d\times {\mathbb R}^d)\setminus (D^c\times D^c)}(u(x)-u(y))(v(x)-v(y)) k(x,y)\, dy \, dx\, ,
\end{equation}
where  $u,v:{\mathbb R}^d\to {\mathbb R}$. This form has recently attracted quite some interest, cf.~\cite{BGPR17, DRV17, FKV15, MSW19, Voi} where different question related to the form were studied. 
In particular, \cite{DRV17} introduces  a new ''non-local normal derivative''
\begin{equation}\label{e:non-local-derivative}
{\mathcal N} u(x)=\int_D (u(x)-u(y))k(x,y)\, dy,\qquad x\in {\mathbb R}^d\setminus \overline{D}\, ,
\end{equation}
with the aim  to solve the Neumann problem
$$
\left\{\begin{array}{rl}
(-\Delta)^{\alpha/2}u=f & \textrm{in } D, \\
{\mathcal N} u=0 & \textrm{in }{\mathbb R}^d\setminus \overline{D},
\end{array}\right.
$$
for the fractional Laplacian $(-\Delta)^{\alpha/2}$, as well as the corresponding heat equation with homogeneous Neumann conditions.  
The authors offer the following probabilistic interpretation of the Neumann heat equation:
\emph{ 
\begin{itemize}
\item[(1)] The solution $u(x,t)$ of the Neumann heat equation is the probability distribution of the position of a particle moving randomly inside $D$.
\item[(2)] When the particle exits $D$, it immediately comes back into $D$.
\item[(3)] The way in which it comes back inside $D$ is the following: If the particle has gone to $x\in {\mathbb R}^d\setminus D$, it may come back to any point $y\in D$, the probability
density of jumping from $x$ to $y$ being proportional to $k(x,y)$.
\end{itemize}
}
In view of the fact that the variational structure of the Neumann problem involves the symmetric bilinear form $\widehat{{\mathcal E}}(u,v)$
where $u$ and $v$ are functions defined on all of ${\mathbb R}^d$ (and not only on $D$), the above probabilistic interpretation is somewhat ambiguous. One goal of this note is to construct two stochastic processes, one living on ${\mathbb R}^d$, the other on $D$, which might fit the intended interpretation. To be more precise, for $x\in D^c$, let $\mu(x):=\int_D k(x,y)dy$, define the measure $m(dx):={\mathbf 1}_D(x)dx +{\mathbf 1}_{D^c}(x)\mu(x)dx$, and set  ${\mathcal F}:=\{u\in L^2({\mathbb R}^d,m(dx)), \widehat{{\mathcal E}}(u,u)<\infty\}$. We will show that $(\widehat{{\mathcal E}}, {\mathcal F})$ is a quasi-regular Dirichlet form on $L^2({\mathbb R}^d,m(dx))$, hence there is a Markov process $Y$ on ${\mathbb R}^d$ (more precisely, on ${\mathbb R}^d\setminus \partial D$) properly associated with $(\widehat{{\mathcal E}}, {\mathcal F})$. The behavior of $Y$ can be described as follows: starting in $D$, the process moves as the isotropic stable process until the first exit time from $D$. At the exit time, it jumps out of $D$ according to the kernel $k(x,y)$. It sits at the exit point $y$ for an exponential time with mean one, then jumps back to $D$ according to probability distribution $k(y,x)/k(y)$ and starts afresh. By deleting the part of this process which lives in $D^c$, we get a process with state space $D$.

The other goal of the note is to look at the corresponding Dirichlet-to-Neumann operator. In the context of the Laplace operator, the classical Dirichlet-to-Neumann operator can be roughly described as follows: take a function $\phi$ defined on the boundary $\partial D$ (for simplicity, here we do not specify the function space that $\phi$ belongs to). Let $u$ be the solution of the Dirichlet boundary value problem (for the Laplacian) with the boundary value $\phi$. Let $\psi$ be the normal derivative of $u$. The mapping $\phi\mapsto \psi$ is called the Dirichlet-to-Neumann operator. There exists a substantial amount of literature on the classical Dirichlet-to-Neumann problem, the results depending on the roughness of the domain and the appropriate function spaces, see for example \cite{AtE11, AtE17, tEO13, tEO17} for the functional-analytic approach and \cite{BV17} for a probabilistic approach. In this note, we solve the Dirichlet-to-Neumann problem for non-local operators, both probabilistically and analytically, and show that the problem is simpler than the one for local operators. In view of the non-locality of the underlying operator (such as the fractional Laplacian), the boundary $\partial D$ is replaced by the exterior $\overline{D}^c$ and the Dirichlet-to-Neumann operator is the mapping between functions defined on $\overline{D}^c$. In the analytic approach,  we define the Dirichlet-to-Neumann operator on $L^2(D^c, \mu(x)dx)$. For the Dirichlet-to-Neumann operator (related to the fractional Laplacian) on different function spaces we refer to \cite{GSU17}. Another closely related research is \cite{FJ01, JS99} where the authors study functions spaces and Dirichlet forms of subordinate reflected diffusions on the closure $\overline{D}$ of a (smooth) open set $D$ in the Euclidean space. They obtain a Weyl decomposition which is the key to construction of a Dirichlet-to-Neumann operator. Their methods of using the spectral synthesis techniques are close in spirit to our approach.

Organization of the paper: in the next section we introduce the singular kernel $k(x,y)$ which generalizes $|x-y|^{-d-\alpha}$ and recall the relevant function spaces based on this kernel. Then we briefly discuss the Dirichlet problem for the corresponding non-local operator and define the Dirichlet-to-Neumann operator in an analytic way.  In Section \ref{s:dirichlet-form} we prove that $(\widehat{{\mathcal E}}, {\mathcal F})$ is a quasi-regular Dirichlet form on $L^2({\mathbb R}^d,m(dx))$, explain the behavior of the corresponding process $Y$ with state space ${\mathbb R}^d$ and argue that its trace on $D^c$ gives a probabilistic interpretation of the Dirichlet-to-Neumann operator. In the last section, starting from $Y$, we construct a stochastic process $Z$ with the state space $D$ and calculate its bilinear form. This last process fits the description from \cite{DRV17} that \emph{when the particle exits $D$, it immediately comes back into $D$.}  Finally, in the appendix, we provide an alternative construction of the process $Y$ and compute its bilinear form.

\section{Preliminaries on function spaces}\label{s:prelim}
In this section we introduce the singular kernel as the L\'evy density of the underlying stochastic process $X$. We also recall several function spaces related to the process   $X$ (or the kernel) that were studied in \cite{FKV15} (see also \cite{DRV17}). Instead of the fractional Laplacian, we will work with a more general non-local operators, in fact generators of isotropic L\'evy processes.

Let $(X_t, {\mathbb P}_x)$ be a L\'evy process on ${\mathbb R}^d$, $d\ge 2$, with the characteristic exponent $\Phi$ of the form
$$
\Phi(\xi)=\int_{{\mathbb R}^d\setminus \{0\}} (1-e^{i\langle \xi, x\rangle}-i\langle \xi, x\rangle \mathbf{1}_{\{|x|< 1\}})\nu(dx),
$$
that is ${\mathbb E}_x[e^{i\langle \xi, X_t-x\rangle}]=e^{-t\Phi(\xi)}$. Here $\nu(dx)$ is the L\'evy measure of $X$, i.e., $\int_{{\mathbb R}^d}(1\wedge |x|^2)\nu(dx)<\infty$. We will assume that $\nu$ has a strictly positive non-increasing radial density (with respect to Lebesgue measure): $\nu(dx)=j(|x|)dx$ with $j:(0,\infty)\to (0, \infty)$ non-increasing. We introduce the \emph{symmetric} kernel $k(x,y):=j(|x-y|)$, $x,y\in {\mathbb R}^d$. 

Let $D\subset {\mathbb R}^d$ be a bounded open set such that the Lebesgue measure of its boundary $\partial D$ is zero. Recall the following function spaces:
\begin{eqnarray*}
W^{k,2}({\mathbb R}^d)&=&\{v:{\mathbb R}^d\to {\mathbb R};\ v\in L^2({\mathbb R}^d), \iint_{{\mathbb R}^d\times {\mathbb R}^d}(v(x)-v(y))^2 k(x,y)\, dy \, dx<\infty\},\\
W^{k,2}_{D}({\mathbb R}^d)&=&\{v\in W^{k,2}({\mathbb R}^d); \ v=0 \text{ a.e.~on }D^c\},\\
V^{k,2}(D)&=&\{v:{\mathbb R}^d\to {\mathbb R};\  v_{|D}\in L^2(D), \iint_{D\times {\mathbb R}^d}(v(x)-v(y))^2 k(x,y)\, dy \, dx<\infty\}.
\end{eqnarray*}
In \cite[Definition 2.1]{FKV15} these spaces were denoted by $H({\mathbb R}^d; k)$, $H_{D}({\mathbb R}^d; k)$ and $V(D; k)$ respectively.

For $u,v\in W^{k,2}({\mathbb R}^d)$ set
$$
{\mathcal E}(u,v):= \frac12 \iint_{{\mathbb R}^d\times {\mathbb R}^d}(u(x)-u(y)(v(x)-v(y)) k(x,y)\, dy \, dx\, .
$$
Then $({\mathcal E}, W^{k,2}({\mathbb R}^d))$ is the regular Dirichlet form corresponding to the $L^2$-semigroup of the process $X$. Moreover,
$$
\|u\|_{W^{k,2}({\mathbb R}^d)}:=\left(\|u\|_{L^2({\mathbb R}^d)}^2+{\mathcal E}(u,u)\right)^{1/2}
$$
is a Hilbert norm on $W^{k,2}({\mathbb R}^d)$. 
 
By symmetry of $k(x,y)$ we have
$$
\int_{D}\int_{D^c}(u(x)-u(y)(v(x)-v(y)) k(x,y)\, dy \, dx = \int_{D^c}\int_{D}(u(x)-u(y)(v(x)-v(y)) k(x,y)\, dy \, dx .
$$
Let
$$
{\mathcal E}^{D}(u,v):= \frac12 \iint_{D\times D}(u(x)-u(y)(v(x)-v(y)) k(x,y)\, dy \, dx\, ,
$$
and
\begin{eqnarray}
\widehat{{\mathcal E}}(u,v)&:=&  \frac12 \iint_{({\mathbb R}^d\times {\mathbb R}^d)\setminus (D^c\times D^c)}(u(x)-u(y)(v(x)-v(y)) k(x,y)\, dy \, dx \label{e:form-whEE}\\
&=&{\mathcal E}^{D}(u,v)+\int_{D}\int_{D^c}(u(x)-u(y)(v(x)-v(y)) k(x,y)\, dy \, dx \, .\nonumber
\end{eqnarray}

In case $k(x,y)=c(d,\alpha)|x-y|^{-d-\alpha}$,  \cite[(3.1), (3.2)]{DRV17} introduces the space $H^s_{D, 0}$ of functions $u:{\mathbb R}^d\to {\mathbb R}$ such that $\|u\|_{L^2(D)}^2+\widehat{{\mathcal E}}(u,u)<\infty$ with the corresponding norm. It is easy to see that $H^s_{D, 0}=V^{k,2}(D)$ in our notation. Denote the corresponding norm by
$$
\|u\|_{V^{k,2}(D)}:=\left (\|u\|_{L^2(D)}^2+\widehat{{\mathcal E}}(u,u)\right)^{1/2}.
$$
Clearly, this is an inner product norm. It is proved in \cite[Proposition 3.1]{DRV17} that $(V^{k,2}(D), \\  \|\cdot \|_{V^{k,2}(D)}) $ is a Hilbert space. Although $k(x,y)=c(d,\alpha)|x-y|^{-d-\alpha}$ in \cite{DRV17}, the proof carries over to $k(x,y)$ as in our setting. Moreover, the proof shows that if $(u_n)_{n\ge 1}$ is a Cauchy sequence in $(V^{k,2}(D), \ \|\cdot \|_{V^{k,2}(D)})$, then a subsequence converges a.e.~in ${\mathbb R}^d$.

Let
$$
V^{k,2}_{D}(D):=\{u\in V^{k,2}(D);\  u=0 \text{ a.e.~in }D^c\}=W^{k,2}_{D}({\mathbb R}^d)\, .
$$
Then $V^{k,2}_{D}(D)$ is a closed subspace of $(V^{k,2}(D), \ \|\cdot \|_{V^{k,2}(D)})$. Indeed, let $(u_n)_{n\ge 1}\subset V^{k,2}_{D}(D)$ and $u=\lim_{n\to \infty}u_n$ in  $(V^{k,2}(D), \ \|\cdot \|_{V^{k,2}(D)})$. Then there is subsequence of $(u_n)$ which converges a.e.~in ${\mathbb R}^d$ to $u$. Hence, $u=0$ a.e.~in $D^c$, i.e.~$u\in V^{k,2}_{D}(D)$.  

Let $u\in L^2_{D}({\mathbb R}^d):=\{v\in L^2({\mathbb R}^d); v=0 \text{ a.e.~on }D^c\}$.  Then
\begin{eqnarray*}
\widehat{{\mathcal E}}(u,u)&=&{\mathcal E}^D(u,u)+\int_D\int_{D^c}(u(x)-u(y))^2k(x,y)\, dy\, dx\\
&\ge &\int_D u(x)^2\left(\int_{D^c}k(x,y)\, dy\right)dx \ge \int_Du(x)^2 \left(\int_{B(x, \mathrm{diam}(D))^c}j(|x-y|)\, dy\right)dx\\
&\ge &c \int_D u(x)^2\, dx\, ,
\end{eqnarray*}
where $c>0$ is a constant not depending on $x\in D$. 
Therefore, for any $u\in L^2_{D}({\mathbb R}^d)$, and in particular for any $u\in V^{k,2}_{D}(D)$, we have that
$$
\|u\|_{L^2(D)}^2\le c^{-1} \widehat{{\mathcal E}}(u,u)<\infty\, .
$$
For $u\in V^{k,2}_{D}(D)$, let $\|u\|_{V^{k,2}_{D}(D)}^2:=\widehat{{\mathcal E}}(u,u)$. The above inequality shows that $(V^{k,2}_{D}(D),\\ \|\cdot\|_{V^{k,2}_{D}(D)})$ is a Hilbert space and the norm $\|\cdot\|_{V^{k,2}_{D}(D)}$ is equivalent to $\|\cdot\|_{V^{k,2}(D)}$.

\begin{lemma}\label{l:functional}
For $u\in V^{k,2}(D)$, let $F:V^{k,2}_{D}(D)\to {\mathbb R}$ be defined by $F(v):=\widehat{{\mathcal E}}(u,v)$, $v\in V^{k,2}_{D}(D)$. 
Then there exists a unique $u_0\in V^{k,2}_{D}(D)$ such that $F(v)=\widehat{{\mathcal E}}(u_0,v)$.
\end{lemma}
\noindent{\bf Proof.} We have that 
$$
|F(v)|=|\widehat{{\mathcal E}}(u,v)|=|{\mathcal E}(u,v)|\le {\mathcal E}(u,u)^{1/2}{\mathcal E}(v,v)^{1/2}={\mathcal E}(u,u)^{1/2}\widehat{{\mathcal E}}(v,v)^{1/2}={\mathcal E}(u,u)^{1/2}\|v\|_{V^{k,2}_{D}(D)}.
$$
This shows that $F$ is a continuous linear functional on the Hilbert space $(V^{k,2}_{D}(D), \|\cdot\|_{V^{k,2}_{D}(D)})$. Hence, there exists a unique $u_0\in V^{k,2}_{D}(D)$ such that $F(v)=\widehat{{\mathcal E}}(u_0,v)$. \qed

Given $\phi :D^c\to {\mathbb R}$, we extend it to all of ${\mathbb R}^d$ by letting $\phi(x)=0$ for $x\in D$. Assume that such extended $\phi\in V^{k,2}(D)$. This is equivalent to 
\begin{equation}\label{e:aux}
\int_{D}\int_{D^c}\phi(x)^2 k(x,y)\, dx\, dy =\int_{D^c}\phi(x)^2 \int_{D} k(x,y)\, dy\, dx<\infty\, .
\end{equation}
Define $\mu:D^c\to [0,\infty]$ by
$$
\mu(x):=\int_{D}k(x,y)\ dy\, ,\qquad x\in D^c\, ,
$$
and note that for $x\in \overline{D}^c$, we have $\mu(x)<\infty$. Indeed, this is clear since for $x\in \overline{D}^c$, $y\in D$, it holds that $|x-y|\ge \mathrm{dist}(x,\partial D)$, hence $\mu(x)\le j(\mathrm{dist}(x,\partial D))|D|$. For $x\in \partial D$, it will usually be the case that $\mu(x)=\infty$.  Let $L^2(D^c, \mu(x)dx)$ be the usual $L^2$ space with the inner product given by
\begin{eqnarray*}
(\phi, \chi)_{L^2(D^c, \mu(x)dx)}&:=&\int_{D^c}\phi(x)\chi(x) \mu(x)\, dx=\int_{\overline{D}^c}\phi(x)\chi(x) \mu(x)\, dx\\
&=&\int_{D^c}\int_{D}\phi(x)\chi(x)k(x,y)\, dy\, dx\, .
\end{eqnarray*}
Extend $\phi\in L^2(D^c, \mu(x)dx)$ to $D$ by letting $\phi(x)=0$, $x\in D$. Then we may regard $L^2(D^c, \mu(x)dx)$ as a (closed) subspace  of $V^{k,2}(D)$ and 
$$
\|\phi\|^2_{L^2(D^c, \mu(x)dx)}=\widehat{{\mathcal E}}(\phi, \phi)\, .
$$

We would like to have a sort of a converse, namely that if $u\in V^{k,2}(D)$, then $u \mathbf{1}_{D^c}\in  V^{k,2}(D)$. Note that
\begin{equation}\label{e:u-restricted}
\widehat{{\mathcal E}}(u \mathbf{1}_{D^c}, u \mathbf{1}_{D^c})=\int_{D^c}\int_{D}(u \mathbf{1}_{D^c}(x)-u \mathbf{1}_{D^c}(y))^2 k(x,y)\, dy\, dx=\int_{D^c}u(x)^2\mu(x)\, dx\, .
\end{equation}
Since the right-hand side need not be finite for $u\in V^{k,2}(D)$, we introduce another function space. Let $m(dx)=\mathbf{1}_{D}(x)dx+\mathbf{1}_{D^c}(x)\mu(x)dx$. Since it may happen that $\int_{D^c}\mu(x)dx=\infty$ (e.g., if $k(x,y)=|x-y|^{d+\alpha}$ for $1<\alpha<2$), the measure $m(dx)$ need not be a Radon measure. 
Set
\begin{eqnarray*}
{\mathcal F}&:=&\{u\in V^{k,2}(D):  \int_{D^c}u(x)^2 \mu(x)\, dx<\infty\}\\
&=&\{u:{\mathbb R}^d\to {\mathbb R}:  u\in L^2({\mathbb R}^d, m(dx)), \widehat{{\mathcal E}}(u,u)<\infty\}\, .
\end{eqnarray*}
For $u\in {\mathcal F}$ define
$$
\|u\|^2_{{\mathcal F}}:=\|u\|_{V^{k,2}(D)}^2+\int_{D^c}u(x)^2\mu(x)\, dx =\|u\|^2_{L^2({\mathbb R}^d, m(dx))}+\widehat{{\mathcal E}}(u,u)\, .
$$
Then $({\mathcal F}, \|\cdot\|_{{\mathcal F}})$ is a Hilbert space. Indeed, the norm is clearly an inner product norm. Suppose that $(u_n)_{n\ge 1}$ is a Cauchy sequence in $({\mathcal F}, \|\cdot\|_{{\mathcal F}})$. Then it is a Cauchy sequence in $(V^{k,2}(D), \|\cdot \|_{V^{k,2}(D)})$. Hence, there exists $u\in V^{k,2}(D)$ such that $\|u_n-u \|_{V^{k,2}(D)}\to 0$. Moreover, a subsequence of $(u_n)$ converges to $u$ a.e.~in ${\mathbb R}^d$, hence also $m$-a.e. Further, $({u_n}_{|D^c})$ is a Cauchy sequence in $L^2(D^c, \mu(x)dx)$, hence converges to some $v\in L^2(D^c, \mu(x)dx)$. Since a subsequence converges to $u$ $m$-a.e., we see that $v=u_{|D^c}$ $m$-a.e. This proves that $u=\lim_{n\to \infty} u_n$ in $({\mathcal F}, \|\cdot\|_{{\mathcal F}})$. Therefore, the following result holds true.

\begin{lemma}\label{l:function-space}
\begin{itemize}
\item [(i)] $({\mathcal F}, \|\cdot\|_{{\mathcal F}})$ is a Hilbert space;
\item[(ii)] If $u\in {\mathcal F}$, then $u\mathbf{1}_{D^c}$ and $u\mathbf{1}_{D}$ are also in ${\mathcal F}$;
\item[(iii)] If $\phi\in L^2(D^c, \mu(x)dx)$, then $\phi$ extended to be zero on $D$ is in ${\mathcal F}$.
\end{itemize}
\end{lemma}
\noindent{\bf Proof.} (ii) If $u\in {\mathcal F}$, then $u\mathbf{1}_{D^c}\in {\mathcal F}$ by \eqref{e:u-restricted}. It follows that $u\mathbf{1}_{D}=u-u\mathbf{1}_{D^c}\in {\mathcal F}$. \qed

\section{Dirichlet-to-Neumann operator}\label{s:d-to-n}
We start this section by recalling the exterior value Dirichlet problem.
Let
$$
{\mathcal L} u(x):= \text{P.V.} \int_{{\mathbb R}^d}(u(x)-u(y))k(x,y)\, dy=\lim_{\epsilon \to 0}\int_{|y-x|>\epsilon}(u(x)-u(y))k(x,y)\, dy\, .
$$
Consider the following exterior value Dirichlet problem (cf.~\cite[Definition 2.5 (D)]{FKV15}). Let $\phi\in L^2(D^c, \mu(x)dx)$. A function $u\in V^{k,2}(D)$ is a solution of 
\begin{equation}\label{e:dp}
\left\{\begin{array}{ll}
	{\mathcal L} u=0 & \text{ in }D\, ,\\
	u=\phi & \text{ on }D^c\, ,
	\end{array}\right.
\end{equation}
if $u-\phi\in W^{k,2}_{D}({\mathbb R}^d)=V^{k,2}_{D}(D)$ and ${\mathcal E}(u, v)=0$ for all $v\in W^{k,2}_{D}({\mathbb R}^d)$. Since for $v\in W^{k,2}_{D}({\mathbb R}^d)$ it holds that $v=0$ a.e.~on $D^c$, this last condition can be written as $\widehat{{\mathcal E}}(u,v)=0$. 

It is shown in \cite[Theorem 3.5 and Theorem 4.4]{FKV15} that there exists a unique solution $u\in V^{k,2}(D)$ of the above Dirichlet problem. Moreover, since $u=\phi$ a.e.~on $D^c$, we see that in fact $u\in {\mathcal F}$. Further, it is also shown in \cite[(3.3)]{FKV15} that there exists $C>0$ such that
$$
\|u\|_{V^{k,2}(D)}\le C \|\phi\|_{L^2(D^c, \mu(x)dx)}\, .
$$
In particular
\begin{equation}\label{e:u-phi}
\widehat{{\mathcal E}}(u,u)^{1/2} \le C\|\phi\|_{L^2(D^c, \mu(x)dx)}\, .
\end{equation}

Let
$$
{\mathcal F}_{D}:=\{u\in {\mathcal F}:\, u=0 \ m-\text{a.e.~in }D^c\},
$$
and note that ${\mathcal F}_{D}=V^{k,2}_D(D)$.
Define now the following simple trace operator. For $u\in {\mathcal F}$ let $\mathrm{Tr}(u):=u_{|D^c}$. Then $\mathrm{ker\,Tr}={\mathcal F}_{D}$. Next we define the space of harmonic functions in $D$ (with respect to the non-local operator ${\mathcal L}$). Let
\begin{eqnarray*}
H(D)&:=&\{u\in {\mathcal F};\ {\mathcal E}(u,v)=0 \text{ for all }v\in \mathrm{ker\  Tr}\}\\
&=&\{u\in {\mathcal F}; \ \widehat{{\mathcal E}}(u,v)=0 \text{ for all }v\in {\mathcal F}_{D}\}.
\end{eqnarray*}
Note that $H(D)$ is a closed subspace of ${\mathcal F}$. This is a consequence of the continuity of the form ${\mathcal E}$, cf.~\cite[Lemma 2.4]{FKV15}. 
We also note that the solution $u\in {\mathcal F}$ of \eqref{e:dp} is in $H(D)$.
Further, ${\mathcal F}_{D}\cap H(D)=\{0\}$. Indeed, if $v$ is in the intersection, then $\widehat{{\mathcal E}}(v,v)=0$. Since $\widehat{{\mathcal E}}(\cdot, \cdot)$ is a norm on ${\mathcal F}_{D}$, it follows that $v=0$.

\begin{lemma}\label{l:decomposition}
 It holds that
$$
{\mathcal F}={\mathcal F}_{D}\oplus H(D)
$$
in the sense that any $u\in {\mathcal F}$ can be uniquely decomposed as $u=v+w$ with $v\in {\mathcal F}_{D}$, $w\in H(D)$ and $\widehat{{\mathcal E}}(v,w)=0$.
Moreover, $\mathrm{ker \, Tr}({\mathcal F})=\mathrm{ker \, Tr}(H(D))$.
\end{lemma}
\noindent{\bf Proof.} By Lemma \ref{l:functional}, there exists a unique $v\in {\mathcal F}_{D}$ such that $\widehat{{\mathcal E}}(u, \psi)=\widehat{{\mathcal E}}(v, \psi)$ for all $\psi\in {\mathcal F}_{D}$. Let $w:=u-v\in {\mathcal F}$. Then for any $\psi \in {\mathcal F}_{D}$ we have $\widehat{{\mathcal E}}(w,\psi)=\widehat{{\mathcal E}}(u,\psi)-\widehat{{\mathcal E}}(v,\psi)=0$. Uniqueness follows from the fact that ${\mathcal F}_{D}\cap H(D)=\{0\}$. The last assertion follows from $\mathrm{ker \, Tr} ({\mathcal F}_{D})=\{0\}$.
\qed

For $u\in {\mathcal F}$, let
$$
{\mathcal N} u(x):=\int_{D}(u(x)-u(y))k(x,y)\, dy\, ,\qquad x\in \overline{D}^c\, ,
$$
cf.~\cite[(1.2)]{DRV17} where (up to  a constant) ${\mathcal N}$ is called a \emph{non-local normal derivative}. Let
$$
\widetilde{{\mathcal N}}u(x):=\frac{{\mathcal N} u(x)}{\mu(x)}=\frac{\int_{D}(u(x)-u(y))k(x,y)\, dy}{\int_{D}k(x,y)\, dy}\, , \qquad x\in \overline{D}^c\, ,
$$
be the normalized non-local normal derivative, see \cite[(3.8)]{DRV17}.
 
We continue by constructing the Dirichlet-to-Neumann operator. 
Let $\phi, \chi:D^c \to {\mathbb R}$ and assume that $\phi, \chi\in L^2(D^c, \mu(x)dx)$.  Let $u,v\in {\mathcal F}$ be the corresponding solutions of the Dirichlet problem \eqref{e:dp}. Since $u-\phi\in W^{k,2}_{D}({\mathbb R}^d)$, we see that $\mathrm{Tr}(u)=\phi$, and similarly, $\mathrm{Tr}(v)=\chi$. Moreover,  $u,v\in H(D)$ .

Define the form ${\mathcal C}:L^2(D^c, \mu(x)dx)\times L^2(D^c, \mu(x)dx)\to {\mathbb R}$ by
$$
{\mathcal C}(\phi, \chi)={\mathcal C}(\mathrm{Tr}(u), \mathrm{Tr}(u)):={\mathcal E}(u,v)=\widehat{{\mathcal E}}(u,v)\, .
$$

By using \eqref{e:u-phi}  in the second inequality below, we see that
$$
|{\mathcal C}(\phi, \chi)|=|\widehat{{\mathcal E}}(u,v)|\le \widehat{{\mathcal E}}(u,u)^{1/2}\widehat{{\mathcal E}}(v,v)^{1/2}\le C^2 \|\phi\|_{L^2(D^c, \mu(x)dx)}\|\chi\|_{L^2(D^c, \mu(x)dx)}\, .
$$
This show that the linear functional $\chi\mapsto {\mathcal C}(\phi, \xi)$ is bounded on $L^2(D^c, \mu(x)dx)$. Hence, there exists $\psi\in L^2(D^c, \mu(x)dx)$ such that ${\mathcal C}(\phi, \chi)=(\psi, \chi)_{L^2(D^c, \mu(x)dx)}$. Define the operator $N:L^2(D^c, \mu(x)dx)\to L^2(D^c, \mu(x)dx)$ by $N\phi=\psi$. Since 
$$
|(N\phi, \chi)|_{L^2(D^c, \mu(x)dx)}=|{\mathcal C}(\phi, \chi)|\le C^2 \|\phi\|_{L^2(D^c, \mu(x)dx)}\|\chi\|_{L^2(D^c, \mu(x)dx)}\, ,
$$
we see that the operator $N$ is bounded. Thus, we have proved the following proposition.

\begin{proposition}\label{p:d-to-n}
There exists a \emph{bounded} operator $N:L^2(D^c, \mu(x)dx)\to L^2(D^c, \mu(x)dx)$ such that
$$
(N\phi, \chi)_{L^2(D^c, \mu(x)dx)}={\mathcal C}(\phi, \chi)=\widehat{{\mathcal E}}(u,v)\, 
$$
for all $\phi, \chi\in L^2(D^c, \mu(x)dx)$.
\end{proposition}

We call $N$ \emph{the Dirichlet-to-Neumann} operator. This is in accordance with the definition of the Dirichlet-to-Neumann operator in the classical setting of the Laplacian, cf.~\cite[p.9]{AtE17}. Another justification is provided by the following observation.  Write $v=v\mathbf{1}_{D^c}+v\mathbf{1}_{D}=v_1+v_2$. Then (since $v_1=0$ on $D$ and $v_1=\chi$ on $D^c$), 
\begin{eqnarray}
\lefteqn{\widehat{{\mathcal E}}(u,v)=\widehat{{\mathcal E}}(u,v_1)} \nonumber \\ 
&=&\int_{D^c}\int_{D}(u(x)-u(y))(v_1(x)-v_1(y))k(x,y)\, dy\, dx \nonumber \\
&=&\int_{D^c} \chi(x) \int_{D}(u(x)-u(y))k(x,y)\, dy\, dx \nonumber \\
&=&\int_{\overline{D}^c} \chi(x) \int_{D}(u(x)-u(y))k(x,y)\, dy\, dx \nonumber \\
&=&\int_{\overline{D}^c} \chi(x) {\mathcal N} u(x)\, dx 
=\int_{D^c} \widetilde{{\mathcal N}}u(x)\chi(x) \mu(x)\, dx =(\widetilde{{\mathcal N}}u, \chi)_{L^2(D^c, \mu(x)dx)}\, .\nonumber
\end{eqnarray}
This shows that $N\phi=\widetilde{{\mathcal N}}u$. 

\section{The Dirichlet form and the corresponding process}\label{s:dirichlet-form}
In this section we show that $(\widehat{{\mathcal E}}, {\mathcal F})$ is a quasi-regular Dirichlet form on $L^2({\mathbb R}^d, m(dx))$ and investigate the corresponding Markov process. Recall that $m(dx)=\mathbf{1}_{D}(x)dx+\mathbf{1}_{D^c}(x)\mu(x)dx$ (and need not be a Radon measure on Borel subsets of ${\mathbb R}^d$), 
$$
{\mathcal F}=\{u:{\mathbb R}^d\to {\mathbb R}:  u\in L^2({\mathbb R}^d, m(dx)), \widehat{{\mathcal E}}(u,u)<\infty\}\, ,
$$
and $\|u\|_{{\mathcal F}}^2=\widehat{{\mathcal E}}(u,u)+\|u\|^2_{L^2({\mathbb R}^d, m(dx))}$. Let us introduce the standard notation in the theory of Dirichlet forms, 
$$
\widehat{{\mathcal E}}_1(u,v):=\widehat{{\mathcal E}}(u,v)+\langle u, v\rangle_{L^2({\mathbb R}^d, m(dx)}\, ,\qquad u,v\in {\mathcal F}\, ,
$$
so that $\|u\|_{{\mathcal F}}^2=\widehat{{\mathcal E}}_1(u.u)$. 
Further note that
$$
{\mathcal F}_{D}=\{u\in {\mathcal F}:\, u=0 \ m-\text{a.e.~in }D^c\}
$$
is equal to $V^{k,2}_{D}(D)$ and $\|\cdot \|_{{\mathcal F}}$ restricted to ${\mathcal F}_{D}$ is equal to $\|\cdot \|_{V^{k,2}(D)}$ restricted to $V^{k,2}_{D}(D)$ (and both are equivalent to $\widehat{{\mathcal E}}(\cdot, \cdot)$). Let $X^{D}$ be the process $X$ killed upon exiting $D$, $(P_t^{D})_{t\ge 0}$ the corresponding $L^2(D, dx)$ semigroup, and $({\mathcal C}^{D}, {\mathcal D}({\mathcal C}^{D}))$ the Dirichlet form. Recall from \cite[Theorem 4.4.3]{FOT} that ${\mathcal D}({\mathcal C}^{D})=\{u\in W^{k,2}({\mathbb R}^d): u=0 \textrm{ a.e.~on }D\}=V_D^{k,2}(D)={\mathcal F}_D$ and that $({\mathcal C}^{D}, {\mathcal F}_D)$ is a regular Dirichlet form. For $u,v\in {\mathcal F}_D$ we have
$$
{\mathcal C}^D(u,v)={\mathcal E}^D(u,v)+\int_D u(x)v(x)\kappa(x)\, dx
$$
where 
$$
\kappa(x)=\int_{D^c}k(x,y)\, dy \, ,\qquad x\in D\,
$$
is the killing function. Also, let
$$
{\mathcal C}^D_1(u,v):={\mathcal C}^D(u,v)+\langle u, v\rangle_{L^2(D,dx)}\, , \qquad u,v \in {\mathcal F}_D\, .
$$

\begin{remark}
{\rm
(i) Let $u\in {\mathcal F}$. Then by Lemma \ref{l:function-space}, $u\mathbf{1}_{D}, u\mathbf{1}_{D^c}\in {\mathcal F}$. Moreover, $u\mathbf{1}_{D}\in {\mathcal F}_D$ while $u\mathbf{1}_{D^c}\in L^2(D^c,\mu(x)dx)$. This show that every $u\in {\mathcal F}$ can be written as a sum of two functions, one from ${\mathcal F}_D$, the other from $L^2(D^c,\mu(x)dx)$. Clearly, such a decomposition is unique. Hence we can write ${\mathcal F}={\mathcal F}_D\oplus L^2(D^c,\mu(x)dx)$. Note that this decomposition is different than the one from Lemma \ref{l:decomposition}.

\noindent
(ii) Note that by Fubini's theorem
$$
\int_{D^c}\mu(y)\, dy=\int_{D^c}\int_D k(x,y)\, dy\, dx=  \int_D \kappa(x)\, dx\, .
$$
This shows that the measure $m(dx)$ is finite if and only if the killing function $\kappa$ is integrable. 
}
\end{remark}
Let $u,v\in {\mathcal F}$ and recall from Lemma \ref{l:function-space} that $u\mathbf{1}_{D}, v\mathbf{1}_{D}, u\mathbf{1}_{D^c}, v\mathbf{1}_{D^c}\in {\mathcal F}$. First note that
$$
\widehat{{\mathcal E}}(u\mathbf{1}_{D}, v\mathbf{1}_{D})={\mathcal E}^D(u,v)+\int_D\int_{D^c}u(x)v(x)k(x,y)\ dy\, dx={\mathcal C}^D(u,v)\, .
$$
We rewrite $\widehat{{\mathcal E}}(u,v)$ now as follows:
\begin{eqnarray*}
\widehat{{\mathcal E}}(u,v)&=&\widehat{{\mathcal E}}(u\mathbf{1}_{D}+u\mathbf{1}_{D^c}, v\mathbf{1}_{D}+v\mathbf{1}_{D^c})\\
&=&\widehat{{\mathcal E}}(u\mathbf{1}_{D}, v\mathbf{1}_{D})+\widehat{{\mathcal E}}(u\mathbf{1}_{D^c}, v\mathbf{1}_{D^c})+\widehat{{\mathcal E}}(u\mathbf{1}_{D}, v\mathbf{1}_{D^c})+\widehat{{\mathcal E}}(u\mathbf{1}_{D^c}, v\mathbf{1}_{D})\\
&=&{\mathcal C}^D(u,v)+\int_{D^c}u(x)v(x)\mu(x)\, dx -\int_D\int_{D^c}\big(u(x)v(y)+u(y)v(x)\big)k(x,y)\, dy\, dx\, .
\end{eqnarray*}
More importantly, we have that
\begin{eqnarray}
\widehat{{\mathcal E}}(u,u)&=&{\mathcal E}^D(u,u)+\int_D\int_{D^c}(u(x)-u(y))^2\, dy\, dx\nonumber\\
&\le &{\mathcal E}^D(u,u)+2\int_D\int_{D^c} (u(x)^2+u(y)^2)\, dy \, dx\nonumber\\
&=&{\mathcal E}^D(u,u) +2\int_D u(x)^2\kappa(x)\, dx+2\int_{D^c}u(y)^2 \mu(y)\, dy \nonumber\\
&\le & 2\left({\mathcal C}^D(u,u)+ \int_{D^c}u(x)^2 \mu(x)\, dx)\right) .\label{e:df-estimate}
\end{eqnarray}

\begin{proposition}\label{p:dirichlet-form}
$(\widehat{{\mathcal E}}, {\mathcal F})$ is a Dirichlet form on $L^2({\mathbb R}^d, m(dx))$.
\end{proposition}
\noindent{\bf Proof.} Clearly, $\widehat{{\mathcal E}}$ is a symmetric bilinear form. Next, we argue that ${\mathcal F}$ is dense in $L^2({\mathbb R}^d, m(dx))$. Let $u\in L^2({\mathbb R}^d, m(dx))$ and write $u=u{\mathbf{1}}_D+u{\mathbf 1}_{D^c}$. The function $u{\mathbf 1}_{D^c}$ is already in ${\mathcal F}$, cf.~Lemma \ref{l:function-space} (iii). Next consider $u{\mathbf{1}}_D\in L^2(D, dx)$. Since ${\mathcal D}({\mathcal C}^{D})={\mathcal F}_D$ is dense in $L^2(D, dx)$, there exists a sequence $(u_n)_{n\ge 1}\subset {\mathcal F}_D$ such that $u_{|D}=\lim_n u_n $ in $L^2(D; dx)$. Extend $u_n$ to all of ${\mathbb R}^d$ by setting $u_n(x)=0$ for $x\in D^c$. Then $u_n+u{\mathbf 1}_{D^c}\in {\mathcal F}$ and converges to $u$ in $L^2({\mathbb R}^d, m(dx))$.

Further, since $({\mathcal F}, \|\cdot\|_{{\mathcal F}})$ is a Hilbert space, the form $(\widehat{{\mathcal E}}, {\mathcal F})$ is closed. Finally, let $v$ be a normal contraction of $u\in {\mathcal F}$. Then $|v(x)-v(y)|\le |u(x)-u(y)|$ for all $x,y\in {\mathbb R}^d$, and thus clearly $\widehat{{\mathcal E}}(v,v)\le \widehat{{\mathcal E}}(u,u)<\infty$, hence normal contraction operates on $\widehat{{\mathcal E}}$. \qed

\begin{remark}
{\rm Note that the Dirichlet form $(\widehat{{\mathcal E}}, {\mathcal F})$ is a special case of \cite[Example 1.2.4]{FOT}. Indeed, let 
$$
j(x,dy)=\frac12 {\mathbf 1}_D(x)k(x,y)dy +\frac12 {\mathbf 1}_{D^c}(x){\mathbf 1}_D(y)\frac{k(x,y)}{\mu(x)}dy\, .
$$
Then $j(x,dy)$ satisfies (j.1), (j.2) and (j.3) from \cite[Example 1.2.4]{FOT}. Further, if $J(dx, dy):=j(x,dy)m(dx)$, then $J(dx,dy)$ is a symmetric measure and 
$$
J(dx,dy)={\mathbf 1}_D(x)k(x,y)dy \, dx +{\mathbf 1}_{D^c}(x){\mathbf 1}_D(y) k(x,y) dy\, dx\, .
$$
Hence
\begin{eqnarray*}
\lefteqn{\iint_{{\mathbb R}^d\times {\mathbb R}^d}(u(x)-u(y))(v(x)-v(y))J(dx, dy)}\\
&=&\frac12 \iint_{({\mathbb R}^d\times {\mathbb R}^d)\setminus (D^c\times D^c)}(u(x)-u(y))(v(x)-v(y))k(x,y)dy\, dx=\widehat{{\mathcal E}}(u,v)\, .
\end{eqnarray*}
Since ${\mathcal F}$ is dense in $L^2({\mathbb R}^d, m(dx))$ (see the proof above), it follows from \cite[Example 1.2.4]{FOT} that $(\widehat{{\mathcal E}}, {\mathcal F})$ is a Dirichlet form.
}
\end{remark}

Next we will show that the form $(\widehat{{\mathcal E}}, {\mathcal F})$ is quasi-regular. For all unexplained notions (such as nests and quasi-continuity) we refer the reader to \cite{CF}.

\begin{theorem}\label{t:qrdf}
The form $(\widehat{{\mathcal E}}, {\mathcal F})$ is a quasi-regular Dirichlet form on $L^2({\mathbb R}^d, m(dx))$. 
\end{theorem}
\proof For a closed subset $F\subset {\mathbb R}^d$, let ${\mathcal F}_F:=\{u\in {\mathcal F}:\, u=0 \text{ $m$-a.e on }{\mathbb R}^d\setminus F\}$. We check that the three properties of \cite[Defintion 1.3.8.]{CF} are satisfied.

\noindent
(i) First we show that there exists an $\widehat{{\mathcal E}}$-nest $(F_j)_{j\ge 1}$ of compact sets. Since $({\mathcal C}^D, {\mathcal F}_D)$ is a regular Dirichlet form, there exists a ${\mathcal C}^D$-nest  $(A_j)_{j\ge 1}$ of compact subset of $D$. This means that $\bigcup_{j\ge 1} {\mathcal F}_{A_j}$ is dense in ${\mathcal F}_D$ with respect to ${\mathcal C}^D_1(\cdot, \cdot)$. Next, let $(B_j)_{j\ge 1}$ be an increasing sequence of compact subsets of $\overline{D}^c$ such that $\bigcup_{j\ge 1}B_j=\overline{D}^c$. For any $u\in L^2(D^c, \mu(x)dx)$, we have that $u{\mathbf 1}_{B_j}\in {\mathcal F}_{B_j}$.  Further, since ${\mathbf 1}_{\overline{D}^c\setminus B_j^c}\rightarrow {\mathbf 1}_{\partial D}=0$ $m$-a.e., by the dominated convergence theorem we get $\|u-u{\mathbf 1}_{B_j}\|_{L^2(D^c, m(dx))}=\|u{\mathbf 1}_{\overline{D}^c\setminus B_j^c}\|_{L^2(D^c, m(dx))}\rightarrow 0$. Set $F_j:=A_j\cup B_j$, $j\ge 1$. Then $F_j$ is a compact subset of ${\mathbb R}^d$. For $u\in {\mathcal F}$ and  $\epsilon>0$ we can find $j\ge 1$, $v\in {\mathcal F}_{A_j}$ and $w\in {\mathcal F}_{B_j}$ such that ${\mathcal C}^D_1(u{\mathbf 1}_D-v,u{\mathbf 1}_D-v)<\epsilon$ and $\langle u{\mathbf 1}_{D^c}-w, u{\mathbf 1}_{D^c}-w\rangle_{L^2(D^c, m(dx))}<\epsilon$. Then $v+w\in {\mathcal F}_{F_j}$ and by \eqref{e:df-estimate}
\begin{eqnarray*}
\lefteqn{\widehat{{\mathcal E}}_1(u-(v+w), u-(v+w))}\\
&=&\widehat{{\mathcal E}}(u-(v+w), u-(v+w))+\langle u-(v+w), u-(v+w)\rangle_{L^2({\mathbb R}^d, m(dx))}\\
&\le &2\left({\mathcal C}^D(u{\mathbf 1}_D-v,u{\mathbf 1}_D-v)+\langle u{\mathbf 1}_{D^c}-w, u{\mathbf 1}_{D^c}-w\rangle_{L^2(D^c, m(dx))}\right)\\
& & +\langle u-(v+w), u-(v+w)\rangle_{L^2({\mathbb R}^d, m(dx))}\\
&\le &4\epsilon\, .
\end{eqnarray*}
This proves that $(F_j)_{j\ge 1}$ is an $\widehat{{\mathcal E}}$-nest of compact sets. Note also that $\bigcap_{j\ge 1}({\mathbb R}^d\setminus F_j)=\partial D$ implying that $\partial D$ is $\widehat{{\mathcal E}}$-polar.

\noindent
(ii) Since $({\mathcal C}^D, {\mathcal F}_D)$ is a regular Dirichlet form, $C_c(D)\cap {\mathcal F}_D$ is ${\mathcal C}^D_1(\cdot, \cdot)$-dense in ${\mathcal F}_D$.   On the other hand, $C_c(\overline{D}^c)$ is dense in $L^2(D^c, \mu(x)dx)$. For $v\in C_c(D)\cap {\mathcal F}_D$, we denote by the same letter the function on ${\mathbb R}^d$ extended to be zero on $D^c$. Similarly, for $w\in C_c(\overline{D}^c)$, the same letter denotes the function on ${\mathbb R}^d$ extended to be zero on $D$. Let ${\mathcal G}:=\{u:\, u=v+w, v\in C_c(D)\cap {\mathcal F}_D, w\in C_c(\overline{D}^c)\}$. Each function in ${\mathcal G}$ is continuous, and therefore $\widehat{{\mathcal E}}$-quasi-continuous. Moreover, ${\mathcal G}$ is dense in $({\mathcal F}, \widehat{{\mathcal E}}_1)$. Indeed, for $u\in {\mathcal F}$ and $\epsilon >0$, there exist $v\in C_c(D)\cap {\mathcal F}_D$ such that ${\mathcal C}^D_1(u{\mathbf 1}_D-v, u{\mathbf 1}_D-v)<\epsilon$, and $w\in C_c(\overline{D}^c)$ such that $\langle u{\mathbf 1}_{D^c}-v, u{\mathbf 1}_{D^c}-v\rangle_{L^2(D^c, m(dx))}<\epsilon$. By \eqref{e:df-estimate}, analogously as in part (i), we get that $\widehat{{\mathcal E}}_1(u-(v+w),u-(v+w))\le 4\epsilon$. 

\noindent
(iii) Let $(A_j)_{j\ge 1}$ be an increasing sequence of compact subsets of $D$ such that $A_j\subset \mathrm{int}(A_{j+1})$ and  $\bigcup_{j\ge 1}A_j=D$. Similarly, let $(B_j)_{j\ge 1 }$ be an increasing sequence of compact subsets of $\overline{D}^c$ such that $B_j\subset \mathrm{int}(B_{j+1})$ and $\bigcup_{j\ge 1}B_j=\overline{D}^c$. Let $v_j\in C_c(D)\cap {\mathcal F}_D$ such that $v_j=1$ on $A_j$ and $v_j=0$ on $D\setminus A_{j+1}$. Similarly, let $w_j\in C_c(\overline{D}^c)$ such that $w_j=1$ on $B_j$ and $w_j=0$ on $D^c\setminus B_{j+1}$. Then $u_j:=v_j+w_j$ is continuous on ${\mathbb R}^d$ and $u_j\in {\mathcal F}$. Thus $(u_k)_{k\ge 1}$ is a family of continuous (hence $\widehat{{\mathcal E}}$-quasi-continuous) functions which clearly separates the points of ${\mathbb R}^d\setminus \partial D$. Since $\partial D$ is $\widehat{{\mathcal E}}$-polar, the third property is verified. \qed

\begin{remark}\label{r:not-regular}
{\rm
It is easy to see that $(\widehat{{\mathcal E}}, {\mathcal F})$ need \emph{not} be a regular Dirichlet form. Indeed, the measure $m(dx)$ need not be a Radon measure which by itself prevents $(\widehat{{\mathcal E}}, {\mathcal F})$ to be regular. Moreover, if the form were regular, then by \cite[Remark 1.3.11.]{CF}, ${\mathcal F}\cap C_c({\mathbb R}^d)$ would separate the points of ${\mathbb R}^d$ (here $C_c({\mathbb R}^d)$ are continuous functions with compact support). But note that if $u\in {\mathcal F}\cap C_c({\mathbb R}^d)$, then also  $u{\mathbf 1}_{D^c}\in {\mathcal F}$, meaning that $\int_{D^c}u(x)^2 \mu(x)\, dx<\infty$. Since $\lim_{x\to \partial D, x\in \overline{D}^c}\mu(x)=+\infty$ (and $u$ is continuous), in case when $m(D^c)=\infty$, this forces $u(x)=0$ for every $x\in \partial D$. Therefore, ${\mathcal F}\cap C_c({\mathbb R}^d)$ does not separate points in $\partial D$. On the other hand, by the general theory, cf.~\cite[Theorem 1.4.3]{CF}, $(\widehat{{\mathcal E}}, {\mathcal F})$ is quasi-homeomorphic to a regular Dirichlet form on a locally compact separable metric space $E$. It is easy to identify the space $E$ and the form: we take $E:={\mathbb R}^d\setminus \partial D=D\cup \overline{D}^c$, the disconnected union of open sets $D$ and $\overline{D}^c$,  and the form is given by the essentially same formula as $\widehat{{\mathcal E}}$:
$$
\widetilde{{\mathcal E}}(u,v)=\iint_{E\setminus (\overline{D}^c\times \overline{D}^c)} (u(x)-u(y))(v(x)-v(y))k(x,y)\, dy \, dx, \qquad u,v\in \widetilde{{\mathcal F}},
$$
where
$$
\widetilde{{\mathcal F}}=\{u:E\to {\mathbb R}:\, \widetilde{{\mathcal E}}(u,u)<\infty\}.
$$
Then $(\widetilde{{\mathcal E}}, \widetilde{{\mathcal F}})$ is a regular Dirichlet form on $L^2(E, m_{|E})$. 
}
\end{remark}

Since  $(\widehat{{\mathcal E}}, {\mathcal F})$ is a quasi-regular Dirichlet form on $L^2({\mathbb R}^d, m(dx))$ we can state the following theorem.
\begin{theorem}\label{t:process-existence}
There exists a Hunt process $Y=(Y_t, {\mathbb Q}_x)$ on ${\mathbb R}^d\setminus \partial D$ properly associated with $(\widetilde{{\mathcal E}}, \widetilde{{\mathcal F}})$.
\end{theorem}
\noindent{\bf Proof.} We deduce from \cite[Theorem 1.5.2.]{CF} that there exists an $\widehat{{\mathcal E}}$-polar set $N\subset {\mathbb R}^d$, and an $m$-symmetric, $m$-tight special Borel standard process $Y=(Y_t, {\mathbb Q}_x)$ on ${\mathbb R}^d\setminus N$ that is properly associated with $(\widehat{{\mathcal E}}, {\mathcal F})$. By inspecting the proof of \cite[Theorem 1.5.2.]{CF} and using Remark \ref{r:not-regular}, we can conclude that $N=\partial D$ and that $Y$ is a Hunt process on ${\mathbb R}^d\setminus \partial D$ properly associated with $(\widetilde{{\mathcal E}}, \widetilde{{\mathcal F}})$. \qed

Let $\zeta$ denote the lifetime of $Y$. 
It is easy to see that the part process of $Y$ on $D$ is precisely $X^D$, the underlying process $X$ killed upon exiting $D$, and that the part process of $Y$ on $\overline{D}^c$ is the process that sits at its starting point for an exponential amount of time (of parameter 1) and then it dies. The behavior of $Y$ is described as follows: starting from $x\in D$, $Y$ moves as the underlying process $X$ until $\tau_D^Y=\inf\{t>0:\, Y_t\notin D\}$, the first exit time from $D$. If $\tau_D^Y<\zeta$, then $Y$ jumps out of $D$ according to the kernel $k(Y_{\tau_D^Y-}, Y_{\tau_D^Y})$ and  $y:=Y_{\tau_D^Y}\in \overline{D}^c$. Then $Y$ sits at $y$ for an exponential amount of time (of parameter 1) and then jumps back to $D$ according to the probability distribution $k(y,x)/\mu(y)$. Once in $D$, the process starts afresh.

We argue now that under a certain weak assumption, the lifetime $\zeta$ of $Y$ is infinite. Indeed, assume that ${\mathbb P}_x(X_{\tau_D^X}\in \partial D)=0$ for every $x\in D$, i.e., when $X$ exits $D$ it does so by jumping into $\overline{D}^c$. Sufficient conditions for this to hold can be found in \cite{Mil75, Szt00, SV08}. Then also $Y_{\tau_D^Y}\in \overline{D}^c$ and $Y$ spends an exponential time (of parameter 1) in $\overline{D}^c$ before coming back. This will be repeated infinitely many time and since the exponential sitting times are independent, the lifetime has to be infinite.

In the appendix we give an alternative construction of the process $Y$ and verify that the corresponding bilinear form is indeed $\widehat{{\mathcal E}}$.

\begin{remark}
{\rm 
It has been proved in \cite[Lemma 2.20]{Voi} (when $k(x,y)=|x-y|^{-d-\alpha}$) that in case of a smooth open set $D$, $(\widehat{{\mathcal E}}, {\mathcal F})$  is a regular Dirichlet form on $L^2({\mathbb R}^d, dx)$. This means that the corresponding Hunt process $\widetilde{Y}$ can start from any point of $\partial D$. Away from the boundary $\partial D$, $\widetilde{Y}$ behaves like a time-changed process $Y$.
}
\end{remark}

In the remaining part of this section we look at the trace process of $Y$ on $\overline{D}^c$ and revisit the Dirichlet-to-Neumann operator. Let 
$$
B_t:=\int_0^t \mathbf{1}_{(Y_s\in \overline{D}^c)}ds\, ,\qquad \qquad \sigma_t:=\inf \{s>0:\, B_s>t \},
$$
and let $W_t:=Y_{\sigma_t}$ be the trace process. This process is a continuous-time Markov chain in $\overline{D}^c$ which sits at the point $x$ an exponential time with mean one, and then jumps to the point $z$ according to the jump distribution $p(x,z)dz$ that we are now going to compute. Let $P_{D}(y,z)$, $y\in D$, $z\in D^c$ be the Poisson kernel of $X$. Then clearly, $P_{D}(y,z)dz$ is the exit distribution of $Y$ from $D$. Let 
$$
T_{D}:=\inf\{t>0:\, Y_t\in D\}\, .
$$
Then
$$
{\mathbb Q}_x(Y_{T_{D}}\in dy)=\frac{k(x,y)}{\mu(x)}\, dy\, ,
$$
and hence
$$
p(x,z)=\int_{D}{\mathbb Q}_x(Y_{T_{D}}\in dy) P_{D}(y,z)=\int_{D}\frac{k(x,y)}{\mu(x)} P_{D}(y,z)\, dy\,  .
$$
In order to identify $p(x,z)$, let $\phi:\overline{D}^c\to [0,\infty)$ be a measurable function. We probabilistically solve the exterior value Dirichlet problem with the exterior data $\phi$. Thus, let
$$
u(y):=\int_{D^c}P_{D}(y,z)\phi(z)\, dz\, ,\quad y\in D\, ,
$$
be the harmonic extension. Then
\begin{eqnarray*}
\psi(x)&:=&\int_{\overline{D}^c}p(x,z)\phi(z)\, dz=\int_{\overline{D}^c}\phi(z)\int_{D}\frac{k(x,y)}{\mu(x)}P_{D}(y,z)\, dz\, dy\\
&=&\int_{D}\frac{k(x,y)}{\mu(x)}\int_{D^c}P_{D}(y,z)\phi(z)\, dz=\int_{D}\frac{k(x,y)}{\mu(x)}u(y)\, dy .
\end{eqnarray*}
Let $Q\phi:=\psi$. Then
\begin{equation}\label{e:i-Q}
(I-Q)\phi(x)=\phi(x)-\int_{D}\frac{k(x,y)}{\mu(x)}u(y)\, dy =\int_{D}(\phi(x)-u(y))\frac{k(x,y)}{\mu(x)}\, dy=\widetilde{{\mathcal N}}u(x)=N\phi(x)\, .
\end{equation}
Thus, $I-Q=N$ is the Dirichlet-to-Neumann operator. In other words, the operator $-N$ is the infinitesimal generator of the trace process. 

Recall that $m(dx)=\mathbf{1}_D(x)dx+\mathbf{1}_{D^c}(x)\mu(x)dx$. Then $m_{|\overline{D}^c}$ is the symmetrizing measure for the kernel $p$: $p(x,dz)m(dx)=p(z,dx)m(dz)$. Indeed, denote the Green function of $D$ with respect to $X$  by $G_{D}$. Then 
\begin{eqnarray*}
p(x,z)\mu(x)&=&\int_{D}k(x,y)P_{D}(y,z)\, dy=\int_{D}k(x,y)\int_{D}G_{D}(y,w)k(w,z)\, dw\, dz\\
&=&\int_{D}k(z,w)\int_{D}G_{D}(w,y)k(y, x)\, dy\, dw =p(z,x)\mu(z)\,  .
\end{eqnarray*}

We end this section by recovering \cite[Corollary 5.2]{BGPR17}. By \eqref{e:i-Q},
\begin{eqnarray*}
\widehat{{\mathcal E}}(u,u)&=&\int_{D^c}\widetilde{{\mathcal N}}u(x)\phi(x)\, dx=\int_{D^c}\left(\phi(x)-\int_{D^c}p(x,z)\phi(z)\, dz\right)\phi(x)\mu(x)\, dx\\
&=&\iint_{D^c\times D^c} \phi(x)^2\mu(x)p(x,z)\, dz \, dx - \iint_{D^c\times D^c} \phi(x)\phi(z)\mu(x)p(x,z)\, dz \, dx\\
&=&\frac12 \iint_{\overline{D}^c\times \overline{D}^c}(\phi(x)-\phi(z))^2 \mu(x)p(x,z)\, dz\, dx.
\end{eqnarray*}
Therefore we have,
\begin{eqnarray*}
{\mathcal E}(u,u)&=&\widehat{{\mathcal E}}(u,u)+\frac{1}{2}\iint_{\overline{D}^c\times \overline{D}^c}(u(x)-u(y))^2k(x,y)\, dy\, dx\\
&=&\frac12 \iint_{\overline{D}^c\times \overline{D}^c} (\phi(x)-\phi(y))^2\left(\mu(x)p(x,y)+k(x,y)\right)\, dy\, dx.
\end{eqnarray*}

\section{The process $Z$}
Now we transform the process $Y$ so that the resulting process lives only in $D$. This new process corresponds to the description in \cite{DRV17} of a process that after it jumps from $D$ immediately returns to $D$.

Let $C_t:=\int_0^t 1_{(Y_s\in D)}\, ds$ be the time $Y$ spends in $D$ until the fixed time $t$. Then $C$ is a positive continuous additive functional whose support is $D$ (cf.~\cite[(5.1.21)]{FOT}). 
Let $\tau_t:=\inf\{s>0:\, C_s>t\}$ be the right-continuous inverse of $C$. Define the new process $Z=(Z_t)_{t\ge 0}$ by $Z_t:=Y_{\tau_t}$. This construction amounts to deleting from the path of $Y$ the part that $Y$ spends in $\overline{D}^c$. The process $Z$ is a right Markov process with the state space $D$ (cf.~\cite[p.175]{CF}).  For a non-negative Borel function $u$ on ${\mathbb R}^d$, let
$$
Hu(x):={\mathbb E}_x[u(Y_{T_{D}}), T_{D}<\infty]=\left\{\begin{array}{ll}
	u(x), & x\in D, \\
	\int_{D}\frac{k(x,y)}{\mu(x)}u(y)\, dy, & x\in \overline{D}^c,\\
	0, &x\in\partial D\, .
	 \end{array}\right.
$$
Further, let 
$$
\left\{ \begin{array}{l}
\widecheck{{\mathcal F}}=\{\phi\in L^2(D, dx):\,  \phi=u \text{ a.e.~on } D \text{ for some }u\in {\mathcal F}_e\}\\
\widecheck{{\mathcal E}}(\phi, \phi)=\widehat{{\mathcal E}}(Hu, Hu):\,  \phi\in \widecheck{{\mathcal F}}, \phi=u \text{ a.e.~on } D, u\in {\mathcal F}_e\, .
\end{array} \right.
$$
Here ${\mathcal F}_e$ denotes the extended Dirichlet space.
By \cite[Theorem 5.2.2. and Theorem 5.2.7]{CF} $(\widecheck{{\mathcal E}},\widecheck{{\mathcal F}})$ is a quasi-regular Dirichlet form and the process $Z$ is properly associated with it. Moreover, if we regard $(\widehat{{\mathcal E}}, \widehat{{\mathcal F}})$ as a regular Dirichlet form $(\widetilde{{\mathcal E}}, \widetilde{{\mathcal F}})$ on $L^2({\mathbb R}^d\setminus \partial D, m_{|{\mathbb R}^d\setminus \partial D})$, then it follows from \cite[Theorem 6.2.1]{FOT}, that $(\widecheck{{\mathcal E}},\widecheck{{\mathcal F}})$ is in fact a regular Dirichlet form on $L^2(D,dx)$.

We compute now $\widecheck{{\mathcal E}}(\phi, \phi)$. For simplicity, let
$$
\widehat{k}(x,y):=\int_{D^c}\frac{k(x,z)k(z,y)}{\mu(z)}\, dz\, ,\quad x,y\in D\, ,
$$
and note that for any $x\in D$,
\begin{eqnarray}
\int_D \widehat{k}(x,y)\, dy&=&\int_D \int_{D^c}\frac{k(x,z)k(z,y)}{\mu(z)}\, dz\, dy=\int_{D^c}\frac{k(x,z)}{\mu(z)}\left(\int_Dk(z,y)\, dy\right)dz \nonumber\\
&=&\int_{D^c}k(x,z)\, dz=\kappa(x)\, .\label{e:k-hat-integral}
\end{eqnarray}
Hence $\widehat{k}$ is an integrable kernel.

We have
\begin{eqnarray*}
\widecheck{{\mathcal E}}(\phi, \phi)&=&\widehat{{\mathcal E}}(Hu, Hu)\\&=&\frac{1}{2}\int_{D}\int_{D}(Hu(x)-Hu(y))^2k(x,y)\, dy\, dx +\int_{D}\int_{D^c}(Hu(x)-Hu(y))^2 k(x,y)\, dy\, dx\\
&=&\frac{1}{2}\int_{D}\int_{D}(u(x)-u(y))^2 k(x,y)\, dy\, dx +\int_{D}u(x)^2\int_{D^c}k(x,y)\, dy\, dx \\
& & -2\int_{D}u(x)\int_{D^c}Hu(y)k(x,y)\, dy\, dx +\int_{D^c}Hu(y)^2 \int_{D}k(x,y)\, dx\, dy\\
&=&\frac{1}{2}\int_{D}\int_{D}(u(x)-u(y))^2 k(x,y)\, dy\, dx +\int_{D}u(x)^2 \kappa(x)\, dx\\
& & -2\int_{D}u(x)\int_{D^c}\left(\int_{D}\frac{k(y,z)}{\mu(y)}u(z)\, dz\right)k(x,y)\, dy\, dx +\int_{D^c}\left(\int_{D}\frac{k(y,z)}{\mu(y)}u(z)\, dz\right)^2 \mu(y)\, dy\\
&=&\frac{1}{2}\int_{D}\int_{D}(u(x)-u(y))^2 k(x,y)\, dy\, dx +\int_{D}u(x)^2 \kappa(x)\, dx\\
& & -2\int_{D}\int_{D}\left(\int_{D^c} \frac{k(x,y)k(y,z)}{\mu(y)}\, dy\right) u(x)u(z)\, dz\, dx \\
& & +\int_{D^c}\frac{1}{\mu(y)}\left(\int_{D}k(y,z)u(z)\, dz\right)\left(\int_{D}k(y,x)u(x)\, dx\right)\\
&=&\frac{1}{2}\int_{D}\int_{D}(u(x)-u(y))^2 k(x,y)\, dy\, dx +\int_{D}u(x)^2 \kappa(x)\, dx\\
& & -2\int_{D}\int_{D}u(x)u(z)\widehat{k}(x,z)\, dz\, dx +\int_{D}\int_{D}u(x)u(z)\widehat{k}(x,z)\, dz\, dx\\
&=&\frac{1}{2}\int_{D}\int_{D}(\phi(x)-\phi(y))^2 k(x,y)\, dy\, dx +\int_{D}\phi(x)^2 \kappa(x)\, dx -\int_{D}\int_{D}\phi(x)\phi(y)\widehat{k}(x,y)\, dy\, dx\\
&=&{\mathcal C}^{D}(\phi, \phi)-\int_{D}\int_{D}\phi(x)\phi(y)\widehat{k}(x,y)\, dy\, dx\, .
\end{eqnarray*}
Note that the calculation above shows that 
$$
\int_{D}\phi(x)^2 \kappa(x)\, dx -\int_{D}\int_{D}\phi(x)\phi(y)\widehat{k}(x,y)\, dy\, dx=\int_{D}\int_{D^c}(Hu(x)-Hu(y))^2 k(x,y)\, dy\, dx\ge 0\, .
$$
Moreover, by use of \eqref{e:k-hat-integral} we have that
$$
\int_{D}\int_{D}\phi(x)^2 \widehat{k}(x,y)\, dy \, dx =\int_{D}\int_{D}\phi(x)^2 \int_{D^c}\left(\frac{k(x,z)k(z,y)}{\mu(z)}\, dz\right)dy\, dx=\int_{D}\phi(x)^2 \kappa(x)\, dx\, ,
$$
implying
$$
\int_{D}\phi(x)^2 \kappa(x)\, dx -\int_{D}\int_{D}\phi(x)\phi(y)\widehat{k}(x,y)\, dy\, dx=\frac{1}{2}\int_{D}\int_{D}(\phi(x)-\phi(y))^2\widehat{k}(x,y)\, dy\, dx\, .
$$
Therefore, we finally have that
$$
\widecheck{{\mathcal E}}(\phi, \phi)=\frac{1}{2}\int_{D}\int_{D}(\phi(x)-\phi(y))^2 (k(x,y)+\widehat{k}(x,y))\, dy\, dx\,  .
$$

\section{Appendix: An alternative construction of the process $Y$}\label{s_appendix}
In this appendix we give an alternative construction of the process $Y$ and compute the bilinear form. The process $Y$ in $D$ behaves like $X$, once it jumps outside $D$, sits at the landing point $x$ for an exponential time, returns to $D$ according to the normalized measure $k(x,y)dy$, and then starts afresh.

Let $X=(X_t, {\mathbb P}_x)$ be the isotropic L\'evy process in ${\mathbb R}^d$ introduced in Section \ref{s:prelim}, and let $\tau=\tau_{D}=\inf\{t>0: X_t\notin D\}$ be the first exit time from $D$. According to \cite[Theorem 10.3, p.305]{Dy65}, the \emph{stopped} process $X^{\tau}$ is a standard process. Since its lifetime is infinite, $X^{\tau}$ is in fact a Hunt process. Define 
$$
A_t:=\int_0^t \mathbf{1}_{(X^{\tau}_s\in \overline{D}^c)}\, ds\, ,\qquad \qquad M_t:=e^{-A_t}\, 
$$
Then $M=(M_t)$ is a continuous strong multiplicative functional, cf.~\cite[III (3.11)]{BG68}. Moreover, if $(R_n)$ is an increasing sequence of stopping times converging to $R$, then $M_{R_n}\to M_R$ a.s.~on $\{R<\infty\}$, cf.~\cite[III (3.14)]{BG68}. Let $\widehat{X}=(\widehat{X}_t, \widehat{{\mathbb P}}_x)$ be the canonical subprocess of $X^{\tau}$ corresponding to $M$, cf.~\cite[III 3.]{BG68}. By \cite[III (3.16 Corollary)]{BG68}, $\widehat{X}$ is a Hunt process. Note that the lifetime $\zeta$ of $\widehat{X}$ is finite almost surely. We see that $\widehat{X}$ is the process that behaves as $X$ while in $D$. If the first exit from $D$ is in $\overline{D}^c$ (which happens a.s.), then $\widehat{X}$ sits at the exit place $X_{\tau}$ for an exponential time with mean one, and then it dies. If the exit place is on $\partial D$ (which will be the case ${\mathbb P}_x$-a.s. for a regular point $x\in \partial D$ ), then $\widehat{X}$ sits at the exit point forever. Also, if $X$ starts at a regular point $x\in \partial D$, then $\widehat{X}$ stays at $x$ forever.

Now we use the piecing out procedure from \cite{INW66}. The \emph{instantaneous distribution} $\mu(\omega, dy)$ is defined by
$$
\mu(\omega, dy):=\frac{\mathbf{1}_{D}(y)k(\widehat{X}_{\zeta-}(\omega), y)dy}{k(\widehat{X}_{\zeta-}(\omega))}\, .
$$
That is, once $\widehat{X}$ is killed it reappears at $y\in  D$ according to the normalized jumping kernel $k(x,y)$ (where $x=\widehat{X}_{\zeta-}$, $y\in D$). Let $Y=(Y_t, {\mathbb Q}_x)$ be the process constructed in \cite[Theorem 1.1]{INW66} by the piecing out procedure. The lifetime of the process $Y$ is infinite. (This is clear because $\widehat{X}$ will be resurrected infinitely many times, and the sum of independent exponential random variables, each with mean 1, is infinite.) Moreover, the lifetime $\zeta$ of $\widehat{X}$ is a totally inaccessible stopping time. Therefore, according to \cite[Corollary of Proposition 4.2]{INW66}, $Y$ is a Hunt process. 

In the remaining part of this appendix we compute the bilinear form of the Hunt process $Y$. For a bounded measurable $u:{\mathbb R}^d\to R$ let $Q_t u(x):={\mathbb Q}_x u(Y_t)$, $t\ge 0$ be the semigroup of $(Y_t, {\mathbb Q}_x)$.  
Recall that $X^{D}$ is the process $X$ killed upon exiting $D$, $(P_t^{D})_{t\ge 0}$ the corresponding $L^2(D, dx)$ semigroup, $({\mathcal C}^{D}, {\mathcal D}({\mathcal C}^{D}))$ the Dirichlet form, and  ${\mathcal D}({\mathcal C}^{D})={\mathcal F}_D$.
The killing function $\kappa:D\to [0,\infty) $ was defined as
$$
\kappa(x):=\int_{D^c}k(x,y)\, dy, \quad x\in D\, .
$$
Then for $u\in  {\mathcal D}({\mathcal C}^{D})={\mathcal F}_D$,
\begin{eqnarray*}
{\mathcal C}^{D}(u,u)&=\frac12 &\int_{D}\int_{D} (u(x)-u(y))^2 k(x,y)\, dy\, dx +\int_{D}\kappa(x)u(x)^2\, dx\\
&=& {\mathcal E}^{D}(u,u)+\int_{D}u(x)^2\left(\int_{D^c}k(x,y)\, dy\right)\, dx\, .
\end{eqnarray*}

\begin{proposition}\label{p:bilinear-Y}
Let $u,v:{\mathbb R}^d\to {\mathbb R}$ be bounded function such that  $u_{|D}\in C(D)\cap {\mathcal F}_D$, $v_{|D}\in {\mathcal F}_D$, $u_{|D^c}\in L^2(D^c,\mu(x)dx)$ and  $v_{|D^c}\in L^1(D^c, \mu(x)dx)$. Then
\begin{equation}\label{e:bilinear-Y}
\lim_{t\to 0}\frac{1}{t}\int_{{\mathbb R}^d}(u(x)-Q_t u(x))v(x)\, m(dx)=\widehat{{\mathcal E}}(u,v)\, .
\end{equation}
\end{proposition}
\noindent{\bf Proof.} 
Clearly, $u\in L^2({\mathbb R}^d, m(dx))$.
Recall that $T_{D}=\inf\{t>0:\, Y_t\in D\}$ and 

\begin{equation}\label{e:distribution-of-YT}
{\mathbb Q}_x(Y_{T_{D}}\in dy)=\frac{k(x,y)}{\mu(x)}\, dy=:n(x,y)\, dy\, .
\end{equation}

Under ${\mathbb Q}_x$, $x\in \overline{D}^c$, $T_D$ has an exponential distribution with mean 1. Thus for $x\in \overline{D}^c$,
\begin{eqnarray*}
{\mathbb Q}_x u(Y_t)&=&{\mathbb Q}_x[u(Y_t), t<T_{D}]+{\mathbb Q}_x[u(Y_t), t\ge T_{D}]\\
&=&u(x){\mathbb Q}_x(t<T_D)+{\mathbb Q}_x[u(Y_t), t\ge T_{D}]\\
&=&e^{-t}u(x)+{\mathbb Q}_x[u(Y_t), t\ge T_{D}]\, ,
\end{eqnarray*}
hence
\begin{equation}\label{e:uxinOmega-c}
u(x)-Q_t u(x)=u(x)(1-e^{-t})-{\mathbb Q}_x[u(Y_t), t\ge T_{D}]\, .
\end{equation}
By \eqref{e:distribution-of-YT} and the fact that $(Y_{T_{D}+t})_{t\ge 0}$ is independent of $T_{D}$, we get
\begin{eqnarray*}
{\mathbb Q}_x[u(Y_t), t\ge T_{D}]&=&\left(\int_{D}\frac{k(x,y)}{\mu(x)}{\mathbb Q}_y(u(Y_t))\, dy\right){\mathbb Q}_x(t<T_{D})\\
&=&(1-e^{-t}) \int_{D}{\mathbb Q}_y(u(Y_t)) n(x,y)\, dy\, .
\end{eqnarray*} 
By the assumption on $u$ and $v$ ($u_{|D}$ bounded and continuous, $v_{|D^c}\in L^1(\mu(x)dx)$), the use of the dominated convergence theorem below is justified, and we get
\begin{eqnarray}\label{e:Omega-c-part2}
\lefteqn{\lim_{t\to 0}\frac{1}{t}\int_{D^c}v(x) {\mathbb Q}_x[u(Y_t), t\ge T_{D}]m(dx)}\nonumber \\
&=&\lim_{t\to 0}\frac{1-e^{-t}}{t}\int_{D^c} v(x)\int_{D}n(x,y){\mathbb Q}_y(u(Y_t))\, dy\, \mu(x)dx\nonumber \\
&=&\int_{D^c}v(x)\int_{D}n(x,y)\mu(x)u(y)\, dy \, dx=\int_{D^c}\int_{D}v(x)u(y)k(x,y)\, dy \, dx \, .
\end{eqnarray}
Further,
\begin{equation}\label{e:Omega-c-part1}
\lim_{t\to 0}\frac{1-e^{-t}}{t}\int_{D^c}u(x)v(x) m(dx)=\int_{D^c}u(x)v(x)\mu(x) \, dx= \int_{D^c}\int_{D}u(x)v(x) k(x,y)\, dy\, dx\, .
\end{equation}
It follows from \eqref{e:uxinOmega-c}-- \eqref{e:Omega-c-part1} that
\begin{eqnarray}\label{e:limit-xinOmega-c}
\lefteqn{\lim_{t\to 0}\frac{1}{t}\int_{D^c}(u(x)-Q_t u(x))v(x)\, m(dx)}\nonumber \\
&=&\int_{D^c}\int_{D}u(x)v(x) k(x,y)\, dy\, dx -\int_{D^c}\int_{D}v(x)u(y)k(x,y)\, dy \, dx \nonumber \\
&=&\int_{D}\int_{D^c}u(y)v(y) k(x,y)\, dy\, dx -\int_{D}\int_{D^c}u(x)v(y)k(x,y)\, dy \, dx
\end{eqnarray}

Now assume that $x\in D$ and let $\tau_D=\inf\{t>0: \, X_t\notin D\}$. Then we have
\begin{eqnarray*}
{\mathbb Q}_x u(Y_t)&=&{\mathbb Q}_x[u(Y_t), t<\tau_{D}]+{\mathbb Q}_x[u(Y_t), t\ge\tau_{D}]\\
&=&{\mathbb E}_x[u(X_t), t<\tau_{D}]+{\mathbb Q}_x[u(Y_t), t\ge\tau_{D}]\\
&=&P_t^{D}u(x)+{\mathbb Q}_x[u(Y_t), t\ge\tau_{D}]\, ,
\end{eqnarray*}
hence
\begin{equation}\label{e:uxinOmega}
u(x)-Q_t u(x)=u(x)-P_t^{D} u(x)-{\mathbb Q}_x[u(Y_t), t\ge\tau_{D}]\, .
\end{equation}

Since $u_{|D}, v_{|D}\in {\mathcal F}_D={\mathcal D}({\mathcal C}^{D})$, we get
\begin{eqnarray}\label{e:limit-xinOmega-1}
\lefteqn{\lim_{t\to 0}\int_{D}(u(x)-P_t^{D}u(x))v(x)\, dx ={\mathcal C}^{D}(u,v)}\nonumber \\
&  =&\frac{1}{2}\int_{D}\int_{D}(u(x)-u(y))(v(x)-v(y) k(x,y)\, dy \, dx \\
& &+ \int_{D}u(x)v(x)\left(\int_{D^c}k(x,y)\, dy\right)\, dx\, .
\end{eqnarray}

Finally, we consider ${\mathbb Q}_x[u(Y_t), t\ge\tau_{D}]$. Let $\mathbf{e}$ be an independent exponential random variable as in the construction of the process $\widehat{X}$ (i.e., $\mathbf{e}$ is the waiting time in $D^c$ before jumping back to $D$). Then
\begin{eqnarray}\label{e:xinOmega-aux}
\lefteqn{{\mathbb Q}_x[u(Y_t), t\ge\tau_{D}]={\mathbb Q}_x[u(Y_t), \tau_{D}\le t <\tau_{D}+\mathbf{e}]+{\mathbb Q}_x[u(Y_t), t\ge\tau_{D}, t\ge\tau_{D}+\mathbf{e}]} \nonumber \\
&=& {\mathbb Q}_x[u(Y_{\tau_{D}}), \tau_{D}\le t <\tau_{D}+\mathbf{e}]+{\mathbb Q}_x[u(Y_t),  t\ge\tau_{D}+\mathbf{e}] \nonumber \\
&=&  {\mathbb Q}_x[u(Y_{\tau_{D}}), \tau_{D}\le t] -{\mathbb Q}_x[u(Y_{\tau_{D}}),t\ge\tau_{D}+\mathbf{e}]+{\mathbb Q}_x[u(Y_t),  t\ge\tau_{D}+\mathbf{e}] \nonumber \\
&=&  {\mathbb E}_x[u(X_{\tau_{D}}), \tau_{D}\le t] -{\mathbb Q}_x[u(Y_{\tau_{D}}),t\ge\tau_{D}+\mathbf{e}]+{\mathbb Q}_x[u(Y_t),  t\ge\tau_{D}+\mathbf{e}] 
\end{eqnarray}

Note that 
$$
{\mathbb Q}_x (\tau_{D}+\mathbf{e}\le t)=\int_0^te^{-s}{\mathbb Q}_x(\tau_{D}+s\le t)\, ds
=e^{-t}\int_0^t e^s\ {\mathbb Q}_x(\tau_{D}\le s)\, ds\,  .
$$
Hence, by right-continuity of $s\mapsto {\mathbb Q}_x(\tau_{D}\le s)$ and the fact that ${\mathbb Q}_x(\tau_{D}\le 0)=0$, we get
\begin{equation}\label{e:limit-exp}
\lim_{t\to 0}\frac{1}{t}{\mathbb Q}_x (\tau_{D}+\mathbf{e}\le t)=\lim_{t\to 0}\frac{1}{t}e^{-t}\int_0^t e^s\ {\mathbb Q}_x(\tau_{D}\le s)\, ds=0\, .
\end{equation}
Since $u$ is bounded,  \eqref{e:limit-exp} implies that
\begin{equation}\label{e:limit-xinOmega-2}
\limsup_{t\ge 0} \frac{1}{t}\big| {\mathbb Q}_x[u(Y_t),  t\ge\tau_{D}+\mathbf{e}]\big| \le \|u\|_{\infty}\lim_{t\to 0}\frac{1}{t}{\mathbb Q}_x (\tau_{D}+\mathbf{e}\le t)=0\, ,
\end{equation}
and, similarly, 
\begin{equation}\label{e:limit-xinOmega-3}
\lim_{t\to 0}\frac{1}{t} {\mathbb Q}_x[u(Y_{\tau_{D}}),  t\ge\tau_{D}+\mathbf{e}]=0\, .
\end{equation}
In order to handle the term ${\mathbb E}_x[u(X_{\tau_{D}}), \tau_{D}\le t] $ we will use the compensation formula
$$
{\mathbb E}_x\sum_{0<s\le t\wedge \tau_D} F(X_{s-}, X_s)={\mathbb E}_x \int_0^{t\wedge \tau_D} \int_{{\mathbb R}^d}F(X_s, y)k(X_s,y)\, dy\, ds
$$
with $F(x,y)=\mathbf{1}_{D}(x)\mathbf{1}_{D^c}(y)u(y)$. Then 
\begin{eqnarray*}
\lefteqn{{\mathbb E}_x[u(X_{\tau_{D}}), \tau_{D}\le t] ={\mathbb E}_x\sum_{0<s\le t\wedge \tau_D} F(X_{s-}, X_s)}\\
&=&{\mathbb E}_x \int_0^{t\wedge \tau_D} \int_{D^c}\mathbf{1}_{(X_s\in D)}u(y)k(X_s,y)\, dy \, ds\\
&=&\int_0^t \int_{D^c} {\mathbb E}_x(\mathbf{1}_{(X_s\in D, s< \tau_D)}k(X_s,y))u(y)\, dy\, ds\, .
\end{eqnarray*}
Since for $x\in D$ and $y\in D^c$, $\lim_{s\to 0}{\mathbb E}_x[\mathbf{1}_{(X_s\in D, s<\tau_D)}k(X_s,y)]=k(x,y)$, we get that
\begin{equation}\label{e:limit-xinOmega-4}
\lim_{t\to 0}\frac{1}{t}{\mathbb E}_x[u(X_{\tau_{D}}), \tau_{D}\le t]=\int_{D^c}u(y)k(x,y)\, dy\, .
\end{equation}
Now it follows from \eqref{e:limit-xinOmega-2}-\eqref{e:limit-xinOmega-4} that
\begin{equation}\label{e:limit-xinOmega-5}
\lim_{t\to 0}\frac{1}{t}\int_{D}{\mathbb Q}_x[u(Y_t), t\ge\tau_{D}] v(x)\, dx = \int_{D}\int_{D^c}v(x)u(y)k(x,y)\, dy\, dx\, .
\end{equation}
Now \eqref{e:uxinOmega}, \eqref{e:limit-xinOmega-1} and \eqref{e:limit-xinOmega-5} imply that
\begin{eqnarray}\label{e:limit-xinOmega}
\lefteqn{\lim_{t\to 0}\frac{1}{t}\int_{D}(u(x)-Q_t u(x))v(x)\, dx=\frac{1}{2}\int_{D}\int_{D}(u(x)-u(y))(v(x)-v(y)) k(x,y)\, dy \, dx} \nonumber \\
& & + \int_{D}\int_{D^c}u(x)v(x) k(x,y)\, dy\, dx -\int_{D}\int_{D^c}v(x)u(y)k(x,y)\, dy\, dx\, .
\end{eqnarray}
Putting together \eqref{e:limit-xinOmega-c} and \eqref{e:limit-xinOmega} we obtain
\begin{eqnarray*}
\lefteqn{\lim_{t\to 0}\frac{1}{t}\int_{{\mathbb R}^d}(u(x)-Q_t u(x))v(x)\, m(dx)}\\
&=&\frac{1}{2}\int_{D}\int_{D}(u(x)-u(y))(v(x)-v(y)) k(x,y)\, dy \, dx\\
& & +\int_{D}\int_{D^c} (u(x)-u(y))(v(x)-v(y)) k(x,y)\, dy \, dx =\widehat{{\mathcal E}}(u,v)\, .
\end{eqnarray*}
\qed

\bigskip
\noindent
{\bf Acknowledgements:} 
I am grateful to Tom ter Elst for several useful discussions about the classical Dirichlet-to-Neumann operator and for the hospitality at the University of Auckland where this project was initiated. Thanks are also due to Renming Song for a helpful input related to the proof of Proposition \ref{p:bilinear-Y}, and Moritz Kassmann and Guy Fabrice Foghem Gounoue for a number of suggestions.
\bigskip
\noindent

\vspace{.1in}
\begin{singlespace}


\end{singlespace}

\vskip 0.1truein

\parindent=0em
\bigskip

{\bf Zoran Vondra\v{c}ek}

Department of Mathematics, Faculty of Science, University of Zagreb, Zagreb, Croatia

Email: \texttt{vondra@math.hr}

\end{document}